\documentclass{amsart}
\usepackage{latexsym,amssymb}
\numberwithin{equation}{section}
\newcommand\Vainerman{Va\u\i nerman}
\newcommand\Szlachanyi{Szlach\'anyi}
\newcommand{\E}{\operatorname E}
 \newcommand{\lsum}{{\textstyle\sum}}
\newtheorem{Lem}{Lemma}[section]
\newtheorem{Prop}[Lem]{Proposition}
\newtheorem{Cor}[Lem]{Corollary}
\newtheorem{Thm}[Lem]{Theorem}

\theoremstyle{definition}

\theoremstyle{remark}
\newtheorem{Rem}[Lem]{Remark}

\renewcommand\o{\otimes}
\newcommand\FS{Frobenius-separable}
\newcommand\IFS{idempotent Frobenius system}
\newcommand\einsch[1]{|_{#1}}
\newcommand\ol{\overline}
\newcommand\tR{\times_R}

\newcommand\nt{\diamond}

\DeclareMathOperator\Hom{\operatorname{Hom}}
\DeclareMathOperator\End{\operatorname{End}}
\DeclareMathOperator\id{\operatorname{id}}
\DeclareMathOperator\Id{\operatorname{\mathit {Id}}}

\newcommand\op{{\operatorname{op}}}

\newcommand\cop{{\operatorname{cop}}}
\newcommand\bop{{\operatorname{bop}}}

\newcommand\LMod[1]{{_{#1}\mathcal M}}
\newcommand\RMod[1]{\mathcal M_{#1}}

\newcommand\Bimod[2]{{_{#1}\mathcal M_{#2}}}

\newcommand\LComod[1]{{^{#1}\mathcal M}}

\newcommand{\ou}[1]{\mathrel{\mathop{\otimes}_{#1}}}

\newcommand{\nti}[1]{\mathrel{\mathop{\diamond}_{#1}}}

\newcommand\sw[1]{{}_{(#1)}}
\newcommand\swm[1]{{}_{(-#1)}}
\newcommand\so[1]{^{(#1)}}

\newcommand\se[1]{{}_{[#1]}}
\newcommand\sem[1]{{}_{[-#1]}}
\newcommand\inv{^{-1}}
\renewcommand\epsilon\varepsilon
\newcommand\noqed{\renewcommand\qed{}}
\newcommand\role{r\^ole}

\makeatletter
\def\namelabel#1#2{\@bsphack
  \protected@write\@auxout{}%
         {\string\newlabel{#1.nme}{{#2}{#2}}}%
  \@esphack}
\def\nmlabel#1#2{\label{#2}\namelabel{#2}{#1}}
\newcommand\nmref[1]{\ref{#1.nme}\ \ref{#1}}
\makeatother

\newlength{\itemwidth}
\newenvironment{imini}{\addtolength\itemwidth{-\leftmargin}\addtolength\itemwidth{-\labelwidth}\\\begin{minipage}{\itemwidth}}{\end{minipage}\\}

\begin{document}
\title{Weak Hopf algebras and quantum groupoids}
\author{Peter Schauenburg}
\address{Mathematisches Institut der Universit\"at M\"unchen, 
Theresienstr.~39, 80333~M\"unchen, Germany}
\email{schauen@rz.mathematik.uni-muenchen.de}

\begin{abstract}
  We give a detailed comparison between the notion of a weak Hopf
  algebra (also called a quantum groupoid by Nikshych and \Vainerman),
  and that of a $\tR$-bialgebra due to Takeuchi (and also called
  a bialgebroid or quantum (semi)groupoid by Lu and Xu). 
  A weak bialgebra is the same
  thing as a $\tR$-bialgebra in which $R$ is \FS. We extend the
  comparison to cover module and comodule theory, duality, and
  the question when a bialgebroid should be called a Hopf
  algebroid.
\end{abstract}
\maketitle
\section{Introduction}
Quantum groupoids (or Hopf algebroids) are to groupoids what 
quantum groups (or Hopf algebras) are to groups:
A Hopf algebroid is the noncommutative analog of the 
function algebra on a groupoid.

A groupoid is a small category, and has a set of morphisms and a 
set of objects (in other terminology arrows and vertices).
Thus the definition of a quantum groupoid should involve two
algebras, one of which (say $H$) plays the \role\ of the function 
algebra
on the quantum space of morphisms, and the other (say $R$)
the \role\ of
the function algebra on the quantum space of objects. Since there is a
source and target assigned to each arrow, one should also expect
(in the reversed direction) two maps from $R$ to $H$ to be part
of the structure, while composition in the groupoid, a partially
defined map on the product, should correspond to a comultiplication
into a suitably defined tensor product of $H$ with itself.

In this note we will compare in detail two notions of quantum 
(semi)groupoids: The $\tR$-bialgebras defined by
 Takeuchi \cite{Tak:GAAA},
and the weak bialgebras defined by B\"ohm and \Szlachanyi\ 
\cite{BohSzl:CCQGND}.
Thus we shall, as it were, provide reference 
[NS] in \cite{HauNil:ITQHA}; while my joint paper with Florian Nill
is announced there optimistically as a preprint to appear shortly, in
reality it was never finished. 
The main result is as follows: 
A weak bialgebra is the same thing as a $\tR$-bialgebra in which
the algebra $R$ is \FS.

The $\tR$-bialgebras defined by Takeuchi \cite{Tak:GAAA}, following
work of Sweedler \cite{Swe:GSA}, are the first quantum (semi)groupoids
appearing in the literature. One should note, though, that
Takeuchi did not consider the analogy with groupoids at all,
whereas this was the key motivation for the definitions of
Lu \cite{Lu:HAQG} and Xu \cite{Xu:QG}, which turn out to be 
equivalent to Takeuchi's, mostly by a translation of notations,
though some care has to be taken about the somewhat different
definitions of counits. See the paper
of Brzezi\'nski and Militaru \cite{BrzMil:BBD} for details.

Weak Hopf algebras were defined by B\"ohm and Szlach\'anyi
\cite{BohSzl:CCQGND}, see also the recent survey 
\cite{NikVai:FQGA} by Nikshych and \Vainerman\
and the literature cited there.
A weak bialgebra $H$ is a coalgebra and algebra such that the
comultiplication is multiplicative, but does not preserve
the unit; dually the multiplication is not counital. These
two requirements are replaced by certain weakened versions.
In this definition there is in the beginning no auxiliary
algebra $R$ playing the role of the function algebra on 
the set of vertices, but rather two anti-isomorphic ``source
and target counital subalgebras'' are constructed from the axioms.

A special case of weak bialgebras, called face algebras, had 
been defined earlier by Hayashi \cite{Hay:FAIGQGT,Hay:BIFA}.
A face algebra turns out to be precisely the special case of
a weak bialgebra in which the, say target, counital subalgebra
is commutative. In \cite{Sch:FATRB} we have shown that a 
face algebra is precisely the special case of a $\tR$-bialgebra
in which the algebra $R$ is commutative and separable. 

It turns out that one can show by essentially the same 
calculations that a weak bialgebra is precisely the special
case of a $\tR$-bialgebra in which $R$ is \FS. A major difference
is that Hayashi's face algebras involve a commutative separable
base algebra by definition, while one has to show that the 
target counital subalgebra of a weak bialgebra is \FS.

The fact that any weak Hopf algebra is a $\tR$-bialgebra (in fact
a Hopf algebroid in the sense of Lu) has meanwhile been shown by
Etingof and Nikshych \cite[Prop.2.3.1]{EtiNik:DQGR1}, who also
show that the target counital subalgebra is separable
(note however that
the formulas between (10) and (11) there seem to claim that the
Frobenius automorphism is always trivial, which is not the case).
This covers a large part of \nmref{weak-to-Tak}.
However, the antipode is used in 
\cite{EtiNik:DQGR1}, while it is not assumed to exist in 
\nmref{weak-to-Tak}. 
On the other hand
the part of the antipode relevant for the proof
(its restrictions to the source
and target counital subalgebras) is present in any weak bialgebra,
even if it does not possess an antipode; this was proved by
Nill \cite{Nil:AWB} along with the fact that the counital 
subalgebras are Frobenius-separable.

After providing some definitions in Sections \ref{sec:tRb} and 
\ref{sec:FS}, we start the real work in \nmref{sec:weakgen} 
by proving some basic facts about
weak bialgebras, notably that the counital subalgebra of a 
weak bialgebra is \FS. As we acknowledged already, this (and all
the facts proved in \nmref{sec:weakgen}) can be found in the 
literature. However, Etingof and Nikshych \cite{EtiNik:DQGR1}
use antipodes, while Nill's paper \cite{Nil:AWB} consistently
uses the assumption that the weak bialgebra in question is 
finite dimensional. The same general assumption
is used in many places in 
\cite{BohNilSzl:WHAIITCS}. Instead of examining the proofs in 
each situation to convince the reader that the extra assumptions
are not necessary, it seemed easier and more useful to develop
the basic facts that we need from scratch.

In \nmref{sec:weakisTak} we prove that any weak bialgebra $H$ is
a $\tR$-bialgebra (which is \cite[Prop.2.3.1]{EtiNik:DQGR1} if
$H$ is a weak Hopf algebra), and conversely, that any 
$\tR$-bialgebra with \FS\ $R$ is a weak bialgebra.

In \nmref{sec:weakHopf} we adress the question when a weak bialgebra
is a weak Hopf algebra. We show in 
\nmref{Hopf-Hopf}
that a weak Hopf algebra can be
characterized as a weak bialgebra $H$ for which a certain canonical 
map $H\ou{H_t}H\rightarrow\Delta(1)(H\o H)$ is a bijection;
this is analogous to a well-known characterization of ordinary Hopf
algebras.
(We should note that
certain identities for antipodes useful for proving one of the 
implications in \nmref{Hopf-Hopf} can be found in 
\cite{BohNilSzl:WHAIITCS}, again under different assumptions).
This also proves that a weak bialgebra is a weak Hopf algebra
if and only if the associated $\tR$-bialgebra is a $\tR$-Hopf
algebra in the sense of the definition we have given in 
\cite{Sch:DDQG}. One should note that this is 
in general rather different
from the definition of a Hopf algebroid by Lu, which involves an
antipodal anti-algebra map and a certain splitting of the 
epimorphism $H\o H\rightarrow H\ou RH$. Our definition by bijectivity
of a canonical map has the advantage of having a canonical characterization
in terms of properties of the module category of $H$.

In \nmref{sec:duality} we show that the correspondence between
weak bialgebras and $\tR$-bialgebras is compatible with taking
duals (in the finite-dimensional case), and with taking the
respective comodule categories.

\section{$\tR$-bialgebras}\nmlabel{Section}{sec:tRb}

In this section we will recall 
the necessary definitions and notations
on $\tR$-bialgebras. For more details we refer to 
\cite{Swe:GSA,Tak:GAAA,Sch:BNRSTHB}. 

Throughout the paper, $k$ denotes a base field.
Modules,  algebras,  unadorned tensor products etc.\ are 
understood to be over $k$ if nothing else is indicated.

Let $R$ be a $k$-algebra. 
We denote the opposite algebra by $\ol R$, we let
 $R\ni r\mapsto \overline r\in \overline R$ denote the obvious 
$k$-algebra antiisomorphism, and abbreviate the enveloping
algebra $R^e:=R\o\ol R$. We write $r\ol s:=r\o\ol s\in R\o\ol R$
for $r,s\in R$.

For $M,N\in\LMod{R^e}$ we let
$$  \int_r {_{\ol r}M}\o {_rN}:=M\o N\big/\langle \ol rm\o n-m\o rn|r\in R,m\in M,n\in N\rangle$$
and we let
$\displaystyle  \int^r {_{\ol r}M}\o {_rN}\subset M\o N$
denote the $k$-submodule consisting of all elements
$\lsum m_i\o n_i\in M\o N$ satisfying 
$\lsum \ol rm_i\o n_i=\lsum m_i\o rn_i$ for all 
$r\in R$. Variations of the $\int_r$ and $\int^r$ notations, which
are due to MacLane, will be used without further notice.
We abbreviate
$\int_r {_{\ol r}M}\o {_rN}=M\nti RN$ for $M,N\in\LMod{R^e}$.

For two $R^e$-bimodules $M$ and $N$ we let
$$M\tR N:=\int^s\int_r {_{\ol r}M_{\ol s}}\o {_rN_s}.$$ 
If $M,N$ are $R^e$-rings, then so is $M\tR N$, 
with multiplication given by
$(\sum m_i\o n_i)(\sum m'_j\o n'_j)=\sum m_im'_j\o n_in'_j$,
and $R^e$-ring structure
$$R^e\ni r\o\ol s\mapsto r\o\ol s\in M\tR N.$$

For $M,N,P\in\Bimod{R^e}{R^e}$ one defines
$$M\tR P\tR N:=\int^{s,u}\int_{r,t}{_{\ol r}M_{\ol s}}
   \o{_{r,\ol t}P_{s,\ol u}}\o{_tN_u}$$
(where $\int^{s,u}:=\int^s\int^u=\int^u\int^s$). There 
are associativity maps
\begin{align*}
  (M\tR P)\tR N &\overset\alpha\rightarrow M\tR P\tR N\\
  M\tR (P\tR N) &\overset{\alpha'}\rightarrow M\tR P\tR N
\end{align*}
given on elements by the obvious formulas (doing nothing), but which 
need not be isomorphisms. If $M,N$ and $P$ are $R^e$-rings, so
is $M\tR N\tR P$, and $\alpha,\alpha'$ are $R^e$-ring maps.

An $R^e$-ring structure on the algebra
$\E=\End(R)$ is given by 
$r\o\ol s\mapsto(t\mapsto rts)$.
We have, for any $M\in\Bimod{R^e}{R^e}$, two $R^e$-bimodule
maps
\begin{xalignat*}2
  \theta :M\tR\End(R)&\rightarrow M;& m\o f&\mapsto \ol{f(1)}m\\
  \theta' :\End(R)\tR M&\rightarrow M;& f\o m&\mapsto f(1)m.
\end{xalignat*}
which are $R^e$-ring homomorphisms if $M$ is an $R^e$-ring.

A $\tR$-bialgebra $L$ is defined to be
an $R^e$-ring
equipped with a comultiplication, a 
map $\Gamma\colon L\rightarrow L\tR L$ of $R^e$-rings
over $R^e$, and a counit, a map
$C :L\rightarrow \E$ of $R^e$-rings,
such that 
\begin{gather}
  \label{tRcoass}\alpha(\Gamma\tR L)\Gamma=\alpha'(L\tR\Gamma)\Gamma\colon L\rightarrow L\tR L\tR L\\
  \label{tRcounit}
  \theta(L\tR C )\Gamma=\id_L=\theta(C \tR L)\Gamma.
\end{gather}
Note that an $R^e$-ring map $\Gamma\colon L\rightarrow L\tR L$ 
induces
a map $\Gamma_0\colon L\rightarrow L\nt L$ in $\LMod{R^e}$, and
an $R^e$-ring map $C :L\rightarrow\E$ 
induces a map $C _0:L\ni \ell\mapsto C (\ell)(1)\in R$
in $\LMod{R^e}$. 
One checks that $\Gamma$ and $C $ fulfill the equations
\eqref{tRcoass} and \eqref{tRcounit} if and only if
$$(\Gamma_0\nt L)\Gamma_0=(L\nt\Gamma_0)\Gamma_0\colon L\rightarrow L\nt L\nt L$$
and $(C _0\nt L)\Gamma_0=\id_L=(L\nt C _0)$ hold.
These mean that $L$, considered as an $R$-$R$-bimodule via the left
$R^e$-module structure, is an $R$-coring.

For $\tR$-bialgebras we will make use of the 
variants $\Gamma(\ell)=:\ell\se 1\o\ell\se 2\in L\tR L$
and 
$$\alpha(\Gamma\tR L)\Gamma(\ell)=:\ell\se 1\o\ell\se 2\o\ell\se 3\in L\tR L\tR L.$$ 
of usual Sweedler notation (reserving $\ell\sw 1\o\ell\sw 2$ for
usual coalgebra structures).

If $L$ is a $\tR$-bialgebra, then
the  tensor product $M\nti RN$
of $M,N\in\LMod L$ can be endowed with an $L$-module structure
by the usual formula
$\ell (m\o n)=\ell\se 1 m\o\ell\se 2 n$.

The suitable definition of comodules over a $\tR$-bialgebra $L$ 
is as follows:
A left $L$-comodule is an $R$-bimodule
  $M$ together with a map 
  $\lambda\colon M\rightarrow L\tR M$ of $R$-bimodules
  such that
$$\alpha'(L\tR\lambda)\lambda=\alpha(\Gamma\tR M)\lambda\colon M\rightarrow L\tR L\tR M$$
and $\theta'(C \tR M)\lambda=\id_M$ hold.
If we denote by $\lambda_0\colon M\rightarrow L\nt M$ the composition
of $\lambda$ with the inclusion of $L\nt M$ into $L\tR M$, then 
coassociativity is equivalent to 
$(L\nt\lambda_0)\lambda_0=(\Gamma_0\nt M)\lambda_0\colon M\rightarrow L\nt L\nt M$
and $(C _0\nt M)\lambda_0=\id_M$.
We will denote by $\LComod L$ the category of left 
$L$-comodules. We will use Sweedler notation in the form
$\lambda(m)=m\sem 1\o m\se 0$ and
$\alpha(\Gamma\tR M)(m)=m\sem 2\o m\sem 1\o m\se 0$
for $L$-comodules.

The category $\LComod L$ of left $L$-comodules over a $\tR$-bialgebra
is monoidal. The tensor product of $M,N\in\LComod L$ is their
tensor product $M\ou RN$ over $R$, equipped with the comodule
structure
  \begin{align*}
    M\ou RN&\rightarrow L\tR(M\ou RN)\\
    m\o n&\mapsto m\sem 1 n\se{-1}\o m\se 0\o n\se 0
  \end{align*}

\section{Frobenius-separable algebras}\nmlabel{Section}{sec:FS}
In this section we compile a few facts and notations on 
\FS\ algebras, that is, Frobenius algebras that are separable so
that the Frobenius system and the separability idempotent 
coincide. All of the material in this section is certainly folklore.

Let $R$ be a $k$-algebra. Recall that $R$ is Frobenius if there is
a Frobenius system $(\phi,e)$ for $R$, which in turn consists by
definition of a $k$-linear map $\phi\colon R\rightarrow k$, and
an element $e=e\so 1\o e\so 2\in R\o R$ such that
$$\forall r\in R\colon r=\phi(re\so 1)e\so 2=e\so 1\phi(e\so 2r).$$
Equivalently, $R$ is finite dimensional, and there is a 
$k$-linear map $\phi\colon R\rightarrow k$ such that
the bilinear form $B_\phi\colon R\times R\rightarrow k$ given by
$B_\phi(x,y)=\phi(xy)$ is nondegenerate.
It follows that $e\in R\o R$ is a Casimir element in the sense that
$(x\o 1)e=e(1\o x)$ in the algebra $R\o R$ for each $x\in R$: by
nondegeneracy of $B_\phi$ it is sufficient to observe
$\phi(yxe\so 1)e\so 2=yx=\phi(ye\so 1)e\so 2x$ 
for all $y\in R$. Recall that the Frobenius automorphism 
$\theta\colon R\rightarrow R$ defined by the Frobenius system
$(\phi,e)$ is by definition the linear map $\theta\colon R\rightarrow R$
with $\phi(xy)=\phi(y\theta(x))$ for all $x,y\in R$. It is an
algebra automorphism. We have
$(1\o x)e=e(\theta(x)\o 1)$ in $R\o R$ for all $x\in R$, by
the calculation
$\phi(ye\so 1\theta(x))e\so 2=\phi(xye\so 1)e\so 2=xy
 =\phi(ye\so 1)xe\so 2$. It is easy to see that this property 
characterizes $\theta$, so that $(\phi,e)$ is a symmetric 
Frobenius system (i.~e.~$B_\phi$ is symmetric) if and only
if $\theta$ is the identity if and only if 
$(\phi,e\so 2\o e\so 1)$ is a Frobenius system if and only if
$e\so 2\o e\so 1=e$.

If $(\phi,e)$ is a Frobenius system, and $t\in R$ is invertible,
then $(\psi,f)$ defined by $\psi(x)=\phi(tx)$ and
$f=(1\o t\inv)e$, is also a Frobenius system by the calculations
$\psi(xf\so 1)f\so 2=\phi(txe\so 1)t\inv e\so 2=t\inv tx=x$ and
$f\so 1\psi(f\so 2x)=e\so 1\phi(tt\inv e\so 2x)=e\so 1\phi(e\so 2x)=x$.

Conversely, if $(\psi,f)$ is another Frobenius system, define
$t:=\psi(e\so 1)e\so 2$. 
Then $\phi(tx)=\psi(e\so 1)\phi(e\so 2x)=\psi(x)$
for all $x\in R$, further
$e=(1\o t)f$ since 
$f\so 1\phi(tf\so 2x)=f\so 1\psi(f\so 2x)=x=e\so 1\phi(e\so 2x)$
for all $x\in R$. 
Finally $t$ is invertible with inverse $\phi(f\so 1)f\so 2$
since $\phi(f\so 1)f\so 2t=\phi(tf\so 1)f\so 2=\psi(f\so 1)f\so 2=1$
and $t\phi(f\so 1)f\so 2=\phi(e\so 1)e\so 2=1$.

Let $(\phi,e)$ be a Frobenius system. Then 
$e$ is a separability itempotent for $R$ if and only if 
$\nabla(e)=1$, in which case we say that $(\phi,e)$ is an \IFS.
If an \IFS\ exists, we will say that $R$ is \FS. 
While every separable algebra $R$ is automatically
symmetric \cite[Expl.(16.58)]{Lam:LMR}, 
it is not necessarily \FS: If $k$
is a field of characteristic $p>0$, then $M_p(k)$ is separable
(hence ---and also obviously--- Frobenius), but not \FS.
However, if $R$ is a commutative separable $k$-algebra, then
$R$ is Frobenius-separable with respect to the trace
functional $\phi\colon R\rightarrow k$. To see this it suffices
to treat the case where $R$ is a field. 
Let $e$ be a separability idempotent, and write 
$e=\sum_{i=1}^n x_i\o y_i$ with $n$ minimal.
Then the elements $x_i$ generate $R$
as a $k$-space (hence they are a basis). 
For take any $x\in R$, put $t:=x_1\inv x$, and 
consider $\phi\in R^*$
with $\varphi(y_i)=\delta_{1,i}$. Then
$$x=(\id\o\varphi)(\sum tx_i\o y_i)=\sum x_i\varphi(y_it).$$
Similarly the $y_i$ form a basis of $R$.
Now let $\phi(r)$ be the trace of multiplication 
by $r$ as an endomorphism of $R$; this defines $\phi\colon R\rightarrow k$,
and we claim that $(\phi,e)$ is an \IFS. Let $(y^i)$ be the dual basis
of $(y_i)$. Then $\phi(r)=\sum y^i(ry_i)$, so that
$$\phi(xe\so 1)e\so 2
  =\sum\phi(xx_i)y_i
  =\sum x^j(xx_jx_i)y_ix
  =\sum x^j(x_i)y_ixx_j
  =\sum y_jxx_j=x$$
follows for $x\in R$. Similarly $e\so 1\phi(e\so 2x)=x$.

If $(\phi,e)$ and $(\psi,f)$ are two \IFS s, then there is
an invertible $t\in R$ with $e\so 1t\inv e\so 2=1$, such that
$\psi(x)=\phi(tx)$ for all $x$ and $f=(1\o t\inv)e$.
\section{Weak Hopf algebras}\nmlabel{Section}{sec:weakgen}
A weak bialgebra $H=(H,\nabla,\Delta)$
is by definition an algebra and coalgebra $H$
such that the comultiplication $\Delta\colon H\rightarrow H\o H$
is multiplicative, and the following four conditions hold
for all $f,g,h\in H:$
\begin{gather}
  \label{lemo}\epsilon(fgh)=\epsilon(fg\sw 1)\epsilon(g\sw 2h),\\
  \label{rimo}\epsilon(fgh)=\epsilon(fg\sw 2)\epsilon(g\sw 1h),\\
  \label{leco}1\sw 1\o 1\sw 2\o 1\sw 3=(\Delta(1)\o 1)(1\o\Delta(1)),\\
  \label{rico}1\sw 1\o 1\sw 2\o 1\sw 3=(1\o\Delta(1))(\Delta(1)\o 1).
\end{gather}
These four conditions weaken the conditions of multiplicativity of
the counit, and comultiplicativity of the unit, which are not
required in a weak bialgebra (whereas $\epsilon(1)=1$ is an easy
consequence of the axioms). Note that by the symmetries of the 
definition, the opposite $H^\op$, coopposite $H^\cop$, and the
opposite and coopposite (or biopposite) $H^\bop$ are weak bialgebras
as well.

We define the source and target counital maps 
$\epsilon_{s,t}\colon H\rightarrow H$
of a weak bialgebra $H$ to be 
\begin{gather*}
  \epsilon_s(h)=1\sw 1\epsilon(h1\sw 2),\\
  \epsilon_t(h)=\epsilon(1\sw 1h)1\sw 2.
\end{gather*}
And denote their images by $H_{s,t}:=\epsilon_{s,t}(H)$; these are
called the source and target counital subalgebras (see below) of $H$.
We note the variants $\epsilon'_{s,t}$ with
\begin{gather*}
  \epsilon_s'(h)=1\sw 1\epsilon(1\sw 2h),\\
  \epsilon'_t(h)=\epsilon(h1\sw 1)1\sw 2.
\end{gather*}
Obviously these are the source and target counital maps for the
weak bialgebra $H^\op$, which means that general statements on them
will follow from general statements on $\epsilon_{s,t}$ mutatis
mutandis. We'll use in the same way that
$\epsilon_s$ is the target counital map of
$H^{\bop}$.

Note $h\sw 1\epsilon_s(h\sw 2)=h\sw 11\sw 1\epsilon(h\sw 21\sw 2)=h$,
so also
\begin{equation}\label{likeantipode}
h\sw 1\epsilon_s(h\sw 2)=\epsilon_t(h\sw 1)h\sw 2=\epsilon'_s(h\sw 2)h\sw 1=h\sw 2\epsilon'_t(h\sw 1)=h
\end{equation}
for all $h\in H$. Moreover
$$\epsilon_t(1\sw 1h)1\sw 2=\epsilon(1'\sw 11\sw 1h)1'\sw 21\sw 2=\epsilon_t(h)$$

We have
$$1\sw 1\o\epsilon_t(1\sw 2)=1\sw 1\o \epsilon(1'\sw 11\sw 2)1'\sw 2
    =1\sw 1\o\epsilon(1\sw 2)1\sw 3=1\sw 1\o 1\sw 2$$
hence 
$$1\sw 1\o 1\sw 2=\epsilon_s(1\sw 1)\o 1\sw 2=\epsilon_s(1\sw 1)\o \epsilon_t(1\sw 2)$$
and the same identities with $\epsilon_{s,t}$ replaced by $\epsilon_{s,t}'$.
In particular $\Delta(1)\in H_s\o H_t$. It also follows that
$\epsilon_{s,t}(x)=\epsilon_{s,t}'(x)=x$ for all $x\in H_{s,t}$, so 
that $\epsilon_t,\epsilon'_t$ are idempotent projectors onto $H_t$.

The calculation
$$\Delta(\epsilon_t(x))=\epsilon(1\sw 1h)1\sw 2\o 1\sw 3=\epsilon(1\sw 1h)1\sw 21'\sw 1\o 1'\sw 2=\epsilon_t(h)1\sw 1\o 1\sw 2$$
for all $h\in H$ shows the first part of
\begin{equation}\label{Deltatarget}
\forall x\in H_t\colon
\Delta(x)=x1\sw 1\o 1\sw 2=1\sw 1x\o 1\sw 2,
\end{equation}
the second is proved similarly, and 
as a corollary we have
\begin{equation}\label{Deltasource}
\forall x\in H_s\colon\Delta(x)=1\sw 1\o x1\sw 2=1\sw 1\o 1\sw 2x.
\end{equation}

For all $g,h\in H$ we have 
\begin{equation}\label{epsepst}
\epsilon(gh)=\epsilon(g1\sw 2)\epsilon(1\sw 1h)=\epsilon(g\epsilon_t(h))
\end{equation}
by \eqref{rimo}, hence
\begin{equation}\label{eps_t-on-prod}
\epsilon_t(gh)=\epsilon(1\sw 1gh)1\sw 2=\epsilon(1\sw 1g\epsilon_t(h))1\sw 2=\epsilon_t(g\epsilon_t(h)),
\end{equation}
and further
$$g\epsilon_t(h)=\epsilon_t(g\sw 1\epsilon_t(h)\sw 1))g\sw 2\epsilon_t(h)\sw 2
   =\epsilon_t(g\sw 1\epsilon_t(h))g\sw 2=\epsilon_t(g\sw 1h)g\sw 2$$
hence 
\begin{equation}\label{H_t-eps_t}
\forall x\in H_t\forall h\in H\colon x\epsilon_t(h)=\epsilon_t(1\sw 1xh)1\sw 2=\epsilon_t(xh).
\end{equation}
In particular $H_t$ is 
multiplicatively closed; it is a subalgebra because also
$\epsilon_t(1)=\epsilon(1'\sw 11)1'\sw 2=\epsilon(1'\sw 1)1'\sw 2=1.$

For $g,h\in H$ we have 
$$\epsilon_t(g)\epsilon_s(h)=\epsilon(1\sw 1g)1\sw 21'\sw 1\epsilon(h1'\sw 2)
  =\epsilon(1\sw 1g)1'\sw 11\sw 2\epsilon(h1'\sw 2)=\epsilon_s(h)\epsilon_t(g),$$
so that the subalgebras $H_s$ and $H_t$ commute element-wise.

\begin{Lem}Let $H$ be a weak bialgebra.
  The target counital map 
  $\epsilon_t$ induces an algebra antiisomorphism
  $H_s\rightarrow H_t$,
  whose inverse is induced by $\epsilon'_s$.
\end{Lem}
\begin{proof}
  To see that $\epsilon_t$ is an an anti-algebra map we compute
  more generally
  $$\epsilon_t(yh)=\epsilon_t(y\epsilon_t(h))=
     \epsilon_t(\epsilon_t(h)y)=\epsilon_t(h)\epsilon_t(y)$$
  for all $y\in H_s$ and $h\in H$, using \eqref{eps_t-on-prod} and
  \eqref{H_t-eps_t}.

  To prove $\epsilon_s'$ induces an inverse isomorphism to 
  the map induced by $\epsilon_t$, we use that  more generally
  \begin{equation}\label{epseps'}
    \forall h\in H\colon\epsilon_t\epsilon_s'(h)=\epsilon_t(h)
  \end{equation}
  by the calculation 
  $$
    \epsilon_t\epsilon_s'(h)
      =\epsilon(1\sw 1\epsilon'_s(h))1\sw 2
      =\epsilon(1\sw 11'\sw 1\epsilon(1'\sw 2h))1\sw 2
      \overset{\eqref{lemo}}=\epsilon(1\sw 1h)1\sw 2=\epsilon_t(h)
  $$
  Applying this to $H^{\bop}$ yields 
  $\epsilon_s'\epsilon_t(h)=\epsilon_s'(h)$ for all $h\in H$, and
  this taken together with \eqref{epseps'} proves the claim.
\end{proof}
\begin{Prop}\nmlabel{Proposition}{targetFS}
  Let $H$ be a weak bialgebra. Then the target counital
  subalgebra $H_t$ is \FS\ with \IFS\ 
  $$(\epsilon\einsch{H_t},(\epsilon_t\o H)\Delta(1))$$
\end{Prop}
\begin{proof}
  The claimed \IFS\ is given more explicitly by
  $$e=\epsilon_t(1\sw 1)\o 1\sw 2=\epsilon(1'\sw 11\sw 1)1'\sw 2\o 1\sw 2.$$
  We have, for all $x\in H_t$:
  $$\epsilon(xe\so 1)e\so 2
     =\epsilon(1'\sw 11\sw 1)\epsilon(x1'\sw 2)1\sw 2
     \overset{\eqref{rico}}=\epsilon(x1\sw 1)1\sw 2=\epsilon_t'(x)=x$$
  and
  $$e\so 1\epsilon(e\so 2x)
     =\epsilon(1'\sw 11\sw 1)1'\sw 2\epsilon(1\sw 2x)
     \overset{\eqref{leco}}=\epsilon(1'\sw 1x)1'\sw 2=\epsilon_t(x)=x$$
  while $\nabla(e)=1$ is quite obvious.
\end{proof}
It follows that 
\begin{equation}\label{H_t-Casimir}
\forall x\in H_t\colon x\epsilon_t(1\sw 1)\o 1\sw 2=\epsilon_t(1\sw 1)\o 1\sw 2x.
\end{equation}

Applying the Lemma to $H^\bop$ yields that $H_s^\op$ is \FS\ with
\IFS\ $(\epsilon,1\sw 1\o\epsilon_s(1\sw 2))$.
In particular
\begin{equation}\label{H_s-Casimir}
\forall y\in H_s\colon 1\sw 1y\o\epsilon_s(1\sw 2)=1\sw 1\o y\epsilon_s(1\sw 2)
\end{equation}

Applying $\epsilon'_s$ (which is an anti-algebra map restricted to
$H_t$) to the first tensor factor of \eqref{H_t-Casimir}, we
obtain
\begin{equation}\label{rechtsrueber}
  \forall x\in H_t\colon 1\sw 1\epsilon'_s(x)\o 1\sw 2=1\sw 1\o 1\sw 2x
\end{equation}
\section{Weak bialgebras are $\tR$-bialgebras}\nmlabel{Section}{sec:weakisTak}
\begin{Thm}\nmlabel{Theorem}{weak-to-Tak}
  Let $(H,\Delta,\epsilon)$ be a weak bialgebra. 
  Put $R:=H_t$.
  Then the
  structure $(H,\Gamma,C)$ 
  of a $\times_R$-bialgebra on $H$ is given as follows:
  The $R^e$-ring structure of $H$ is given by
  $\iota(x\o\ol y)=x\epsilon'_s(y)$, the comultiplication
  $$\Gamma\colon H\rightarrow H\tR H\subset H\nt H$$
  is the composition of $\Delta$ with the canonical surjection
  $H\o H\rightarrow H\nt H$.
  The counit is
  $$C\colon H\ni h\mapsto (x\mapsto \epsilon_t(hx))\in\End(R).$$
\end{Thm}
\begin{proof}
  $H$ is an $R^e$-ring as claimed since $\epsilon_s'$ induces
  an antiisomorphism of $H_t$ with $H_s$, and $H_s$ and $H_t$ commute
  element-wise.
  
  That $\Gamma_0\colon H\rightarrow H\o H\rightarrow H\nt H$
  takes values in $H\tR H$ follows from 
  $$\Gamma(h)=\Gamma(h\cdot 1)
    =h\sw 11\sw 1\o h\sw 21\sw 2$$
  and \eqref{rechtsrueber}.

  It is clear that $\Gamma$ is an algebra map, since $\Delta$ is
  multiplicative and $\Gamma(1)=1$ in $H\nt H$.
  Also, $\Gamma$ is a map of $R^e$-rings by \eqref{Deltasource}
  and \eqref{Deltatarget},
  and obviously coassociative since $\Delta$ is.

  The map $C$ is unit-preserving since $\epsilon_t$ is idempotent, and
  multiplicative since 
  $$C(g)C(h)(x)=C(g)(\epsilon_t(hx))=\epsilon_t(g\epsilon_t(hx))=\epsilon_t(ghx)$$
  for all $g,h\in H$ and $x\in H_t$, using \eqref{eps_t-on-prod}.
  Moreover
  $C(y)(x)=\epsilon_t(yx)=yx$
   and
  $C(\ol y)(x)=\epsilon_t(\epsilon_s'(y)x)=\epsilon_t(x\epsilon'_s(y))=x\epsilon_t\epsilon'_s(y)=xy$
  for $x,y\in H_t$ show that $C$ is a map of $R^e$-rings.
  
  It remains to check that $C$ is a counit:
  We have 
  $$
    C(h\se 1)(1)h\se 2
    =\epsilon(1\sw 1h\se 1)1\sw 2h\se 2
    =\epsilon(1\sw 1h\sw 1)1\sw 2h\sw 2
    =h
  $$
  as well as
  \begin{multline*}
    \ol{C(h\se 2)(1)}h\se 1
    =\ol{\epsilon(1\sw 1h\se 2)1\sw 2}h\se 1
    =\epsilon(1\sw 1h\sw 2)\epsilon'_s(1\sw 2)h\sw 1
    \\=\epsilon(1\sw 1h\sw 2)1'\sw 1\epsilon(1'\sw 21\sw 2)h\sw 1
    \overset{\eqref{lemo}}=\epsilon(1'\sw 2h\sw 2)1'\sw 1h\sw 1
    =h
  \end{multline*}
  for all $h\in H$.
\end{proof}
The theorem above (which is \cite[Prop.2.3.1]{EtiNik:DQGR1} in the
case where $H$ is a weak Hopf algebra) shows that any weak
bialgebra is a $\tR$-bialgebra in which, by \nmref{targetFS},
 $R$ is \FS. We will also prove a converse \nmref{Tak-to-weak}.
Just as in the case of commutative separable $R$ treated
in  \cite{Sch:FATRB}, this is based
on the following simple observation:
\begin{Rem}\nmlabel{Remark}{keyrem}
  Let $R$ be a separable algebra with separability idempotent
  $e$. Then for $M\in\RMod R$ and $N\in\LMod R$ the
  identity on $M\o N$ induces an isomorphism
  $$\gamma\colon Me\so 1\o e\so 2N\rightarrow M\ou RN$$
  with inverse given by $\gamma\inv(m\o n)=me\so 1\o e\so 2n$.
\end{Rem}
Before using this (implicitly) to prove \nmref{Tak-to-weak}, we
will use it to compare the tensor product defined on modules
over a weak bialgebra by B\"ohm and Szlach\'anyi 
\cite{BohSzl:WHAIIRTDMT}
with the tensor product
defined on the modules over the corresponding $\tR$-bialgebra.
The tensor product on $H$-modules for a weak bialgebra $H$ is
given by $M\odot N:=\Delta(1)(M\o N)$ for $M,N\in\LMod H$, 
with the diagonal left $H$-module structure induced via $\Delta$.
\begin{Prop}\nmlabel{Proposition}{modcateq}
  Let $H$ be a weak bialgebra. Then the isomorphisms
  $$\gamma=\gamma_{MN}\colon M\odot N\rightarrow M\nt N$$
  for $M,N\in\LMod H$ endow the identity functor with 
  the structure of a monoidal functor
  $$(\Id,\gamma)\colon (\LMod H,\nt)\rightarrow(\LMod H,\odot)$$
\end{Prop}
\begin{proof}
The \IFS\ we have found for $R=H_t$ in \nmref{targetFS} is such
that $\ol{e\so 1}\o e\so 2=\Delta(1)$. Thus
$\gamma$ is a vector space isomorphism by \nmref{keyrem};
it is
linear by definition of comultiplication in the $\tR$-bialgebra
associated to the weak bialgebra $H$. Coherence of the monoidal
functor is evident since $\gamma$ is induced by the identity
(and we skip treating unit objects altogether).
\end{proof}
\begin{Rem}
  The arguments used in \nmref{modcateq} could be rewritten 
  to be a different proof of \nmref{weak-to-Tak}: A weak 
  bialgebra $H$ is an $R^e$-ring for $R=R_t$; since $R$ is
  separable, we can use \nmref{keyrem} to endow the underlying
  functor $\LMod H\rightarrow \LMod{R^e}$ with the structure
  of a monoidal functor. It then follows from \cite[Thm.5.1]{Sch:BNRSTHB}
  that $H$ has a $\tR$-bialgebra structure.
\end{Rem}
We now proceed to prove the converse of \nmref{weak-to-Tak}:
\begin{Thm}\nmlabel{Theorem}{Tak-to-weak}
  Let $R$ be a \FS\ algebra with \IFS\ $(\phi,e)$.
  Let $(H,\Gamma,C)$ be a 
  $\times_R$-bialgebra. Then the structure $(H,\Delta,\epsilon)$
  of a 
  weak bialgebra on $H$ is given by
  \begin{gather*}
    \Delta(h)=\sum \ol{e\so 1}h\se 1\o e\so 2h\se 2 \\
    \epsilon(h)=\phi(C(h)(1))
  \end{gather*}
\end{Thm}
\begin{proof}
  The map $\Delta$ is well-defined since 
  $$f\colon H\nt H\ni g\o h\mapsto \ol{e\so 1}g\o e\so 2h\in H\o H$$
  is well-defined, since 
  $\ol{e\so 1}\ol xg\o e\so 2h
   =\ol{xe\so 1}g\o e\so 2h=\ol{e\so 1}g\o e\so 2xh$
  holds for all $g,h\in H$ and $x\in R$. We have
  \begin{multline*}
    \Delta(h\sw 1)\o h\sw 2
      =\ol{e\so 1}(\ol{\tilde e\so 1}h\se 1)\se 1
         \o e\so 2(\ol{\tilde e\so 1}h\se 1)\se 2
         \o \tilde e\so 2h\se 2
      \\=\ol{e\so 1}h\se 1\se 1\o e\so 2\ol{\tilde e\so 1}h\se 1\se 2
         \o\tilde e\so 2h\se 2
      =\ol{e\so 1}h\se 1\o e\so 2\ol{\tilde e\so 1}h\se 2\se 1
         \o\tilde e\so 2h\se 2\se 2
      \\=\ol{e\so 1}h\se 1\o \ol{\tilde e\so 1}
         (e\so 2h\se 2)\se 1\o \tilde e\so 2(e\so 2h\se 2)\se 2
      =h\sw 1\o\Delta(h\sw 2)
  \end{multline*}
  showing that $\Delta$ is coassociative. The map $\epsilon$ is a 
  counit since 
  \begin{multline*}
    h\sw 1\epsilon(h\sw 2)
      =\ol{e\so 1}h\se 1\phi(C(e\so 2h\se 2)(1))
      =\ol{e\so 1}h\se 1\phi(e\so 2C(h\se 2)(1))
      \\=\ol{C(h\se 2)(1)}h\se 1=h
  \end{multline*}
  and
  \begin{multline*}
    \epsilon(h\sw 1)h\sw 2
      =\phi(C(\ol{e\so 1}h\se 1)(1))e\so 2h\se 2
      =\phi(C(h\se 1)(1)e\so 1)e\so 2h\se 2
      \\=C(h\se 1)(1)h\se 2=h
  \end{multline*}

  $\Delta$ is multiplicative by the calculation
  \begin{multline*}
    \Delta(g)\Delta(h)
      =\ol{e\so 1}g\se 1\ol{\tilde e\so 1}h\se 1
            \o e\so 2g\se 2\tilde e\so 2h\se 2
      =\ol{e\so 1}g\se 1h\se 1\o e\so 2\tilde e\so 1\tilde e\so 2h\se 2
      \\=\ol{e\so 1}g\se 1h\se 1\o e\so 2g\se 2h\se 2
      =\Delta(gh)
  \end{multline*}
  for all $g,h\in H$, using $\Gamma(g)\in H\tR H$.

  We have
  \begin{multline*}
    \epsilon(g1\sw 1)\epsilon(1\sw 2h)
      =\epsilon(g\ol{e\so 1})\epsilon(e\so 2h)
      =\phi(C(g\ol{e\so 1})(1))\phi(C(e\so 2h)(1))
      \\=\phi(C(g)(e\so 1))\phi(e\so 2C(h)(1))
      =\phi(C(g)(C(h)(1))=\phi(C(gh)(1))=\epsilon(gh)
  \end{multline*}
  and 
  \begin{multline*}
    \epsilon(g1\sw 2)\epsilon(1\sw 1h)
    =\phi(C(ge\so 2)(1))\phi(C(\ol{e\so 1}h)(1))
    =\phi(C(g)(e\so 2))\phi(C(h)(1)e\so 1)
    \\=\phi(C(g)(C(h)(1)))=\epsilon(gh)
  \end{multline*} 
  for $g,h\in H$, 
  \begin{multline*}
    (H\o \Delta)\Delta(1)
    =\ol {e\so 1}\o\Delta(e\so 2)
    =\ol{e\so 1}\o\ol{\tilde e\so 1}(e\so 2)\se 1\o\tilde e\so 2(e\so 2)\se 2
    \\=\ol{e\so 1}\o\ol{\tilde e\so 1}e\so 2\o\tilde e\so 2
    =(\Delta(1)\o 1)(1\o\Delta(1))
    \\=\ol{e\so 1}\o e\so 2\ol{\tilde e\so 1}\o \tilde e\so 2
    =(1\o\Delta(1))(\Delta(1)\o 1).\qed
  \end{multline*}     
\noqed
\end{proof}

\begin{Rem}
  Let $(H,\Gamma,C)$ be a $\tR$-bialgebra. Then for any \IFS\
  $(\phi,e)$ we obtain a weak bialgebra structure $(H,\Delta_\phi,e_\phi)$
  from \nmref{Tak-to-weak}. 

  On the other hand, if a weak bialgebra structure $(H,\Delta,\epsilon)$
  is given, we obtain a $\tR$-bialgebra structure from 
  \nmref{Tak-to-weak}, along with an \IFS\ for the target counital
  subalgebra $R:=H_t$ from \nmref{targetFS}.

  Assume we start with an \IFS\ on $R$ and a 
  $\tR$-bialgebra $(H,\Gamma,C)$. 
  Consider the weak bialgebra
  $(H,\Delta,\epsilon)$ obtained from it. 
  Assuming that the maps from 
  $R$ and from $\ol R$ to $H$ making $H$ an $R^e$-ring 
  are injective, it is easy to see that $H_t\cong R$, and that
  the \IFS\ on $H_t$ obtained from \nmref{targetFS} is the same 
  as the \IFS\ on $R$ originally given.

  On the other hand, assume we start with a weak bialgebra 
  $(H,\Delta,\epsilon)$,
  and consider the \FS\ algebra $R=H_t$ with \IFS\ $(\phi,e)$
  as in \nmref{targetFS}, and the $\tR$-bialgebra $(H,\Gamma,\epsilon)$
  as in \nmref{weak-to-Tak}. Then for any choice of an \IFS\
  $(\psi,f)$ on $R$ we obtain a weak bialgebra structure
  $(H,\Delta_\psi,\epsilon_\psi)$ from \nmref{Tak-to-weak}. It is 
  quite obvious that $\Delta_\phi=\Delta$ and $\epsilon_\phi=\epsilon$,
  that is, we get the original weak bialgebra back provided we 
  choose the \IFS\ it defines. What happens if we choose another one?
  Then there
  is an invertible $t\in R$ with
  $e\so 1t\inv e\so 2=1$, $\psi(x)=\phi(t x)$ for all $x\in R$, and
  $f=(1\o t\inv)e$, and we obtain
  $$\Delta_\psi(h)=\ol{f\so 1}h\sw 1\o f\so 2h\sw 2=h\sw 1\o t\inv h\sw 2$$
  and 
  $$\epsilon_\psi(h)=\epsilon(tC(h)(1))=\epsilon(t\epsilon_t(h)).
  =\epsilon(th)$$
  This kind of twisting of a weak bialgebra structure by an invertible
  element in the target counital subalgebra is considered by
  Nikshych \cite{Nik:SWHA}. We see that \nmref{weak-to-Tak}
  and \nmref{Tak-to-weak} relate Takeuchi's $\tR$-bialgebras to
  weak bialgebras up to such twists, which corresponds  well to 
  the viewpoint in \cite{Nik:SWHA} that twistings by invertible
  elements in the counital subalgebra should be considered as 
  particularly irrelevant for the structure of $H$. Weak bialgebras
  that are such twists of each other can simply be obtained as 
  different weak bialgebra versions of the same $\tR$-bialgebra.
\end{Rem}

\section{Weak Hopf algebras are $\tR$-Hopf algebras}\nmlabel{Section}{sec:weakHopf}

Etingof and Nikshych have shown that a weak Hopf algebra is a 
Hopf algebroid in the sense of Lu. 

In this section we compare
the weak Hopf algebra axioms to a different notion of 
``Hopf algebroid'', namely that of a $\tR$-Hopf algebra introduced
in \cite{Sch:DDQG}. By definition
\cite[Def.3.5]{Sch:DDQG},
a $\tR$-bialgebra is
a $\tR$-Hopf algebra if and only if the canonical map
$$H\ou{\ol R}H\ni g\o h\mapsto g\se 1\o g\se 2h\in H\nt H$$
is a bijection. This is analogous to a well-known characterization
of ordinary bialgebras. Moreover, the definition is backed in 
\cite{Sch:DDQG} by a characterization of $\tR$-Hopf algebras through
a canonical property of their module categories.

By definition, a weak bialgebra $H$ is a weak Hopf algebra if there
is an endomorphism $S$ of the $k$-space $H$ such that
for all $h\in H$
\begin{gather*}
  S(h\sw 1)h\sw 2=\epsilon_s(h),\\
  h\sw 1S(h\sw 2)=\epsilon_t(h),\\
  S(h\sw 1)h\sw 2S(h\sw 3)=S(h).
\end{gather*}
The axioms imply immediately that 
$$S(h\sw 1)\epsilon_t(h\sw 2)=S(h)=\epsilon_s(h\sw 1)S(h\sw 2).$$
Hence we have, for $x\in H_s$,
\begin{multline*}
\epsilon_s(xh\sw 1)S(h\sw 2)
  =S(h\sw 1)xh\sw 2S(h\sw 3)=S(h\sw 1)x\epsilon_t(h\sw 2)
  \\=S(h\sw 1)\epsilon_t(h\sw 2)x=S(h)x.
\end{multline*}
The antipode is an algebra antihomomorphism by
\begin{multline*}
S(gh)=S(g\sw 1h\sw 1)\epsilon_t(g\sw 2h\sw 2)
       =S(g\sw 1h\sw 1)\epsilon_t(g\sw 2\epsilon_t(h\sw 2))
       \\=S(g\sw 1h\sw 1)g\sw 2\epsilon_t(h\sw 2)S(g\sw 3)
       =S(g\sw 1h\sw 1)g\sw 2h\sw 2S(h\sw 3)S(g\sw 3)
       \\=\epsilon_s(g\sw 1h\sw 1)S(h\sw 2)S(g\sw 2)
       =\epsilon_s(\epsilon_s(g\sw 1)h\sw 1)
        S(h\sw 2)S(g\sw 2)
       \\=S(h)\epsilon_s(g\sw 1)S(g\sw 2)=S(h)S(g)
\end{multline*}
and
$$S(1)=S(1\sw 1)1\sw 2S(1\sw 3)=S(1\sw 1)1\sw 21'\sw 1S(1'\sw 2)
   =\epsilon_s(1)\epsilon_t(1')=1.$$

\begin{Thm}\nmlabel{Theorem}{Hopf-Hopf}
  Let $H$ be a weak bialgebra. Then $H$ is a weak Hopf algebra
  if and only if the map
  $$\beta_0\colon H\o H\ni g\o h\mapsto g\sw 1\o g\sw 2h\in H\o H$$
  induces an isomorphism 
  $$\beta\colon H\ou{H_s} H
    \rightarrow \Delta(1)(H\o H).$$
\end{Thm}
\begin{proof}
  First, assume that $H$ has an antipode $S$.
  Define $\ol\beta_0\colon H\o H\rightarrow H\ou{H_s} H$ by
  $\ol\beta(g\o h)=g\sw 1\o S(g\sw 2)h$. Then
  \begin{multline*}
  \beta\ol\beta_0(g\o h)
   =\beta(g\sw 1\o S(g\sw 2)h)
   =g\sw 1\o g\sw 2S(g\sw 3)h
   =g\sw 1\o\epsilon_t(g\sw 2)h
   \\=1\sw 1g\sw 1\o\epsilon_t(1\sw 2g\sw 2)h
   =1\sw 1g\sw 1\o\epsilon(1'\sw 11\sw 2g\sw 2)1'\sw 2h
   \\=1\sw 1g\sw 1\o\epsilon(1\sw 2g\sw 2)1\sw 3h
   =1\sw 1g\o 1\sw 2h
  \end{multline*}
  and
  \begin{multline*}
    \ol\beta_0\beta_0(g\o h)=g\sw 1\o S(g\sw 2)g\sw 3h
     =g\sw 1\o\epsilon_s(g\sw 2)h
     =g\sw 1\epsilon_s(g\sw 2)\o h
     =g\o h
  \end{multline*}
Thus the restriction of $\ol\beta_0$ is an inverse to $\beta$.

Now assume that $\beta$ has an inverse $\beta\inv$.
Define $\pi\colon H\ou{H_s}H\rightarrow H$ by
$\pi(g\o h)=\epsilon_s(g)h$, and define $S\colon H\rightarrow H$
by $S(h)=\pi\beta\inv(1\sw 1h\o 1\sw 2)$ for $h\in H$.
We claim that $S$ is an antipode for $H$.
For this we first compute
\begin{multline*}
  S(h\sw 1)h\sw 2
    =\pi(\beta\inv(1\sw 1h\sw 1\o 1\sw 2))h\sw 2
    =\pi(\beta\inv(1\sw 1h\sw 1\o 1\sw 2)(1\o h\sw 2))
    \\=\pi(\beta\inv(1\sw 1h\sw 1\o 1\sw 2h\sw 2)
    =\pi\beta\inv(h\sw 1\o h\sw 2)
    =\pi(h\o 1)=\epsilon_s(h).
\end{multline*}
Next, we claim that the inverse of $\beta$ is the restriction of
the map
$$\gamma\colon H\o H\ni g\o h\mapsto g\sw 1\o S(g\sw 2)h\in H\ou{H_s}H.$$
This is verified by the calculation
\begin{multline*}
  \gamma\beta_0(g\o h)
    =\gamma(g\sw 1\o g\sw 2h)
    =g\sw 1\o S(g\sw 2)g\sw 3h
    =g\sw 1\o\epsilon_s(g\sw 2)h
    \\=g\sw 1\epsilon_s(g\sw 2)\o h
    =g\o h.
\end{multline*}
Using, for $y\in H_s$,
\begin{multline*}
  S(yh)=\pi\beta\inv(1\sw 1yh\o 1\sw 2)
    =\pi\beta\inv(1\sw 1h\o 1\sw 2\epsilon_t(y))
    \\=\pi\beta\inv(1\sw 1h\o 1\sw 2)\epsilon_t(y)
    =S(h)\epsilon_t(y),
\end{multline*}
we find
\begin{multline*}
  1\sw 1h\o 1\sw 2
    =\beta\beta\inv(1\sw 1h\o 1\sw 2)
    =\beta(h\sw 1\o S(1\sw 1h\sw 2)1\sw 2)
    \\=\beta(h\sw 1\o S(h\sw 2)\epsilon_t(1\sw 1)1\sw 2)
    =\beta(h\sw 1\o S(h\sw 2))
    =h\sw 1\o h\sw 2S(h\sw 3).
\end{multline*}
and can apply $\epsilon\o H$ to the result to obtain
$\epsilon_t(h)=h\sw 1S(h\sw 2)$. We finish the proof by
calculating
$$S(h\sw 1)h\sw 2S(h\sw 3)=S(h\sw 1)\epsilon_t(h\sw 2)
  =S(h\sw 1)\epsilon_t\epsilon_s'(h\sw 2)
  =S(\epsilon'_s(h\sw 2)h\sw 1)=S(h)$$
for all $h\in H$.
\end{proof}
\begin{Cor}
  Let $H$ be a weak bialgebra. Then the following are equivalent:
  \begin{enumerate}
    \item $H$ is a weak Hopf algebra.
    \item The associated $\tR$-bialgebra $H$ is a 
      $\tR$-Hopf algebra.
  \end{enumerate}
\end{Cor}
\begin{proof}
  The identity induces an isomorphism 
  $\gamma\Delta(1)(H\o H)\rightarrow H\nt H$ by 
  \nmref{modcateq}. The composition $\gamma\beta$ is the map
  $$H\ou{H_s}H\ni g\o h\mapsto g\se 1\o g\se 2h\in H\nt H$$
  required to be bijective in the definition of a $\tR$-Hopf
  algebra.
\end{proof}

For ordinary Hopf algebras, a well-known application of
the characterization \nmref{Hopf-Hopf} is due to Nichols
\cite{Nic:QHA}: Any finite dimensional quotient bialgebra $H/I$
of a Hopf algebra $H$ is itself a Hopf algebra. Dually, every
finite-dimensional subbialgebra of a Hopf algebra is itself
a Hopf algebra. 
Our results will not be quite as striking. We cannot prove
that a finite-dimensional weak subbialgebra $B\subset H$ of a weak
Hopf algebra $H$ is necessarily a weak Hopf algebra. But at least
we can give a criterion purely in terms of the module structure
of $B$ over the source and target counital subalgebras.

To prepare, we note an observation of Nikshych and \Vainerman\
\cite[2.1.12]{NikVai:AVFDQG}:
\begin{Lem}
  Let $f\colon B\rightarrow H$ be a homomorphism of
  weak bialgebras. Then $f$ induces isomorphisms
  $B_t\cong H_t$ and $B_s\cong H_s$
\end{Lem}
\begin{proof}
  We only treat the target counital subalgebra. It is trivial
  to check that $f(B_t)\subset H_t$. We denote the induced
  map $B_t\rightarrow H_t$ by $f$ again. Define 
  $$g\colon H_t\ni x\mapsto \epsilon(xf(1\sw 1))1\sw 2\in B_t.$$
  Then $$gf(x)=\epsilon(f(x)f(1\sw 1))1\sw 2
    =\epsilon(f(x1\sw 1))1\sw 2=\epsilon(x1\sw 1)1\sw 2=x$$
  for all $x\in B_t$, and
  $$fg(x)=\epsilon(xf(1\sw 1))f(1\sw 2)=\epsilon(xf(1)\sw 1)f(1)\sw 2=\epsilon(x1\sw 1)1\sw 2=x$$
  for all $x\in H_t$, so that $g$ is inverse to $f$.
\end{proof}

\begin{Thm}
  Let $H$ be a weak Hopf algebra, and $B\subset H$ a 
  finite-dimensional weak subbialgebra
  (i.e.\ subalgebra and subcoalgebra) of $H$. 
  The following
  are equivalent:
  \begin{enumerate}
    \item $B$ is a weak Hopf algebra.
    \item The right $B_s$-module $B$ is isomorphic to 
      the $B_s$-module $B$ obtained by restricting the left 
      $B_t$-module $B$ along $\epsilon_t$.
  \end{enumerate}
\end{Thm}
\begin{proof}
  As a special case of the preceding Lemma we have 
  $B_t=H_t$ and $B_s=H_s$. 
  For the implication (1)$\Rightarrow$(2) we need to use that
  the antipode of a finite-dimensional quasi-Hopf algebra is 
  bijective \cite[2.10]{BohNilSzl:WHAIITCS}. Now 
  $S$ is an algebra antiautomorphism, and for $y\in H_s$
  $$S(y)=S(y\sw 1)\epsilon_t(y\sw 2)=S(1\sw 1)\epsilon_t(y1\sw 2)
    =S(1\sw 1)1\sw 2\epsilon_t(y)=\epsilon_t(y),$$ 
  so that (2) follows.

  Now assume (2), and fix an isomorphism $f\colon B\rightarrow B$
  satisfying $f(b\epsilon_s'(x))=xf(b)$ for all $b\in B$
  and $x\in H_t$. Then
  $$B\ou{H_s}B\ni b\o c\mapsto c\o f(b)\in B\nt B$$
  is an isomorphism of vector spaces, hence its domain and codomain
  have the same dimension. The canonical map
  $$B\ou{H_s}B\ni b\o c\mapsto b\se 1\o b\se 2c\in B\nt B$$
  is the restriction (note all the modules that occur are
  projective) of the canonical map for $H$, hence injective,
  hence bijective, so that $B$ is a quasi-Hopf algebra.
\end{proof} 
With essentially the same proof we can show:
\begin{Thm}
  Let $H$ be a weak Hopf algebra, and $B=H/I$ a finite-dimensional
  quotient weak bialgebra (i.e.\ $I$ is a coideal and an ideal).
  The following
  are equivalent:
  \begin{enumerate}
    \item $B$ is a weak Hopf algebra.
    \item The right $B_s$-module $B$ is isomorphic to 
      the $B_s$-module $B$ obtained by restricting the left 
      $B_t$-module $B$ along $\epsilon_t$.
  \end{enumerate}
\end{Thm}
\section{Duality}\nmlabel{Section}{sec:duality}
  In \cite[Sec.5]{Sch:DDQG} we have discussed a notion of 
  skew pairing and duality suitable for $\tR$-bialgebras.
  Let $R$ be a $k$-algebra and $H,\Lambda$ two $\tR$-bialgebras.
  We have defined \cite[Def.5.1]{Sch:DDQG}
  a skew pairing between $\Lambda$ and $H$ to be a 
  $k$-linear map
  $\tau\colon\Lambda\o H\rightarrow R$ satisfying
  \begin{alignat}
    \tau((r\o \ol s)\xi(t\o\ol u)|h)v
      &=r\tau(\xi|(t\o\ol v)h(u\o\ol s)),
      \label{skp.1}\\
    \tau(\xi|gh)
       &=\tau(\ol{\tau(\xi\se 2|h)}\xi\se 1|g),
       \label{skp.2}
       &\tau(\xi|1)=C(\xi)(1),\\
    \tau(\xi\zeta|g)
       &=\tau(\xi|\tau(\zeta|g\se 1)g\se 2),
       \label{skp.3}
       &\tau(1|h)=C(h)(1)
  \end{alignat}
  for all $r,s,t,u,v\in R$, $\xi,\zeta\in \Lambda$ and $g,h\in L$.

  As pointed out in \cite{Sch:DDQG}, it is essential that we define
  a skew pairing rather than a pairing in this situation. 
  (An alternative chosen by Kadison and \Szlachanyi\ \cite{KadSzl:DBDTRE}
  is to consider pairings between ``left'' and ``right'' bialgebroids.)

  For 
  weak bialgebras it is no problem to define a Hopf algebra pairing,
  of course, though the problem for $\tR$-bialgebras has its 
  counterpart
  in the fact that the source counital subalgebra of
  the dual of a finite-dimensional weak bialgebra $H$ is 
  canonically isomorphic
  to the target rather than the source counital subalgebra of $H$.

  We define a skew pairing between weak bialgebras $\Lambda,H$ to
  be a linear map $\tau_0\colon \Lambda\o H\rightarrow k$
  satisfying
  \begin{xalignat*}2
    \tau_0(\xi|gh)&=\tau_0(\xi\sw 1|g)\tau_0(\xi\sw 2|h),
    &\tau_0(\xi|1)=\epsilon(\xi),\\
    \tau_0(\xi\zeta|h)&=\tau_0(\zeta|h\sw 1)\tau_0(\xi|h\sw 2),
    &\tau_0(1|h)=\epsilon(h)
  \end{xalignat*}
  for $\xi,\zeta\in\Lambda$ and $g,h\in H$.

  For the rest of the section let $R$ be a \FS\ algebra with \IFS\ 
  $(\phi,e)$.
\begin{Lem}
  Let $\Lambda,H$ be two $\tR$-bialgebras.
 
  If $\tau\colon \Lambda\o H\rightarrow R$ is a skew pairing of
  $\tR$-bialgebras. Then $\tau_0:=\phi\tau\colon\Lambda\o H\rightarrow k$
  is a skew pairing between the corresponding weak bialgebras.
\end{Lem}
\begin{proof}
  We first note that $\ol{e\so 1}h\sw 1\o e\so 2h\sw 2=h\sw 1\o h\sw 2$
  holds in $H\o H$ for all $h\in H$ (and similar formulas for $\Lambda$),
  so that
  $$    \tau(\zeta|h\sw 1)h\sw 2
      =\phi(\tau(\zeta|h\sw 1)e\so 1)e\so 2h\sw 2
      =\phi(\tau(\zeta|\ol{e\so 1}h\sw 1))e\so 2h\sw 2
      =\tau_0(\zeta|h\sw 1)h\sw 2
  $$
  for all $\zeta\in\Lambda$ and $h\in H$.
  It follows that
  \begin{multline*}
    \tau_0(\xi\zeta|h)
      =\phi\tau(\xi\zeta|h)
      =\phi(\tau(\xi|\tau(\zeta|h\se 1)h\se 2))
      \\=\tau_0(\xi|\tau(\zeta|h\sw 1)h\sw 2))
      =\tau_0(\xi|\tau_0(\zeta|h\sw 1)h\sw 2)
  \end{multline*}
  for all $\xi,\zeta\in\Lambda$ and $h\in H$. Moreover
  $$  \tau_0(1|h)=\phi\tau(1|h)=\phi(C(h)(1))=\epsilon(h).$$
  On the other side, using that
  $$\ol{\tau(\xi\sw 2|h)}\xi\sw 1
    =\ol{e\so 1}\phi(e\so 2\tau(\xi\sw 2|h))\xi\sw 1
    =\ol{e\so 1}\xi\sw 1\phi(\tau(e\so 2\xi\sw 2|h))
    =\xi\sw 1\tau_0(\xi\sw 2|h)$$
  for $\xi\in\Lambda$ and $h\in H$, we find
  \begin{multline*}
    \tau_0(\xi|gh)
      =\phi\tau(\xi|gh)
      =\phi\tau(\ol{\tau(\xi\se 2|h)}\xi\se 1|g)
      \\=\tau_0(\ol{\tau(\xi\sw 2|h)}\xi\sw 1|g)
      =\tau_0(\xi\sw 1\tau_0(\xi\sw 2|h)|g)
  \end{multline*}
  for all $\xi\in\Lambda,g,h\in H$, while
  $$\tau_0(\xi|1)=\phi\tau(\xi|1)=\phi(C(\xi)(1))=\epsilon(\xi).$$
\end{proof}
By \eqref{skp.1}, a skew pairing between $\tR$-bialgebras
$\Lambda$ and $H$ defines a map $\Lambda\rightarrow\Hom_{\ol R-}(H,R)$.
In the case that $H$ is finitely generated projective, there
is a unique $\tR$-bialgebra structure on 
$H^\vee:=\Hom_{\ol R-}(H,R)$ such that evaluation defines a 
skew pairing $H^\vee\o H\rightarrow k$, see \cite[Thm.5.13]{Sch:DDQG}.

In our situation, where $R$ is \FS, any $R$-module is projective, so
$H$ is a finitely generated projective left $\ol R$-module if and only
if $H$ is finite dimensional over $k$. 
Then the vector space dual $H^*$ of the weak bialgebra $H$ 
has a natural weak bialgebra structure, consisting of the usual
dual algebra of the coalgebra $H$, and the dual coalgebra of the 
algebra $H$. Note that $(H^*)^\op$ has a skew pairing with $H$.
\begin{Prop}\nmlabel{Proposition}{dual}
  Let $H$ be a finite dimensional $\tR$-bialgebra. Then the 
  weak bialgebra corresponding to $H^\vee$ is the opposite 
  $(H^*)^\op$ of
  the dual $H^*$ of the weak bialgebra corresponding to $H$.
\end{Prop}
\begin{proof}
  Evaluation of $H^\vee$ on $H$ defines a skew pairing of 
  $\tR$-bialgebras which is nondegenerate in its right argument
  by definition of $H^\vee$, and also in its left argument since
  $H$ is finitely generated projective as a left $\ol R$-module.

  We have seen that $\tau$ induces a skew pairing 
  $\tau_0\colon H^\vee\o H\rightarrow k$ of weak bialgebras.
  It only remains to verify that $\tau_0$ is nondegenerate.

  So let first $\xi\in H^\vee$, and assume that
  $\tau_0(\xi|h)=0$ for all $h\in H$. Then for all $h\in H$ we have
  $0=\phi(\tau(\xi|\ol{e\so 1}h))e\so 2
    =\phi(\tau(\xi|h)e\so 1)e\so 2=\tau(\xi|h)$
  and hence $\xi=0$ by definition of $H^\vee$.

  By a parallel argument we can show that $\tau_0$ is nondegenerate
  in the left argument as well.
\end{proof}

\begin{Rem}
  Let $(H,\Delta)$ be a finite-dimensional weak bialgebra with counital
  subalgebra $R$. Let $(H,\Gamma)$ be the associated $\tR$-bialgebra.
  To distinguish, let $\LComod{(H,\Delta)}$ denote the category
  of left comodules over the ordinary $k$-coalgebra $H$, and
  let $\LComod{(H,\Gamma)}$ denote the category of left comodules
  over the $\tR$-bialgebra $H$. By
  \cite[Cor.5.15]{Sch:DDQG} one has an equivalence of monoidal 
  categories $\LComod H\cong\LMod{H^\vee}$. By \nmref{modcateq}
  and \nmref{dual}
  we may replace $\LMod{H^\vee}$ by the module category
  $\LMod{(H^*)^\op}$ over the opposite of the dual weak bialgebra,
  which in turn is the comodule category $\LComod{(H,\Delta)}$ 
  over the weak bialgebra $H$.
  Combining, we have a category equivalence $\LComod{(H,\Delta)}\cong\LComod{(H,\Gamma)}$,
  which we will now derive more directly, and without using finiteness.
\end{Rem}
\begin{Prop}
  Let $(H,\Delta)$ be a weak bialgebra, and $(H,\Gamma)$ the 
  associated $\tR$-bialgebra. 
  \begin{enumerate}
    \item Let $M$ be a comodule over the $\tR$-bialgebra $H$, with
      comodule structure $\lambda\colon M\ni m\mapsto m\sem 1\o m\se 0\in H\tR M$.
      Then the underlying vector space of $M$ is a left $H$-comodule
      over the $k$-coalgebra $H$ with comodule structure
      $$\delta\colon M\ni m\mapsto \ol{e\so 1}m\se 0\o e\so 2m\se 0\in H\o M.$$
    \item Let $M$ be a left $H$-comodule over the $k$-coalgebra $H$,
      with comodule structure map $\delta$.
      Then $M$ is an $R$-bimodule by
      \begin{imini}
      \begin{equation}\label{wiebi}
        rms:=\epsilon(rm\swm 1s)m\sw 0
      \end{equation}
      \end{imini}
      and a left $H$-comodule for the $\tR$-bialgebra $H$ with
      the comodule structure $\lambda$ such that
      $$M\xrightarrow\lambda H\tR M\subset H\nt M$$
      is the composition
      $$M\xrightarrow\delta H\o M\rightarrow H\nt M$$
      in which the second map is the canonical epi.
  \end{enumerate}
  The two constructions describe a bijection between 
  the two types of comodule structures on a given $k$-vector space
  $M$. In particular, one has a category equivalence
  $\LComod{(H,\Delta)}\cong\LComod{(H,\Gamma)}$.
\end{Prop}
\begin{proof}
  If we assume that $M$ is a comodule over the $\tR$-bialgebra 
  $H$, then the calculation showing that it is a comodule over
  the coalgebra $H$ as claimed in (1) is a spitting image of 
  the proof that $H$ is an ordinary coalgebra in 
  \nmref{Tak-to-weak}. 

  So assume that $M$ is a comodule over the
  ordinary coalgebra $H$. 
  The bimodule structure \eqref{wiebi} was first defined by Nill
  \cite[Prop.4.1]{Nil:AWB} under the assumption that $H$ is finite
  dimensional. As a first indication that the 
  structure is appropriate, note
  that we have
  \begin{equation} 
     \label{passtscho}
     \epsilon(rh\sw 1s)h\sw 2=rhs
  \end{equation}
  for all $r,s\in R$ and $h\in H$
  by \eqref{Deltatarget}. To see that \eqref{wiebi} defines a 
  bimodule structure, we only make a sample calculation, say of
  associativity of the left module structure:
  \begin{multline*}
    r(sm)=\epsilon(sm\swm 1)rm\sw 0
         =\epsilon(sm\swm 2)\epsilon(rm\swm 1)m\sw 0
         \\=\epsilon(r\epsilon(sm\swm 1\sw 1)m\swm 1\sw 2)m\sw 0
      \overset{\text{\eqref{passtscho}}}=\epsilon(r(sm\swm 1))m\sw 0
      =(rs)m.
  \end{multline*}
  The calculation for associativity of the right $R$-module structure,
  and compatibility of the left and right module structures, 
  are analogous.
  
  Once we note now (compare \cite[(4.24)]{Nil:AWB}) that
  \begin{multline*}
    m\swm 11\sw 1\o m\sw 01\sw 2
      =m\swm 21\sw 1\o \epsilon(m\swm 11\sw 2)m\sw 0
      \\=(m\swm 1)\sw 1\epsilon((m\swm 1)\sw 2)\o m\sw 0
      =m\swm 1\o m\sw 0
  \end{multline*}
  holds for $m\in M$, the rest of the proof of (2) is again the same
  as the proof that $H$ is a $\tR$-bialgebra in \nmref{weak-to-Tak}.
 
  We omit showing that the two constructions 
  described in (1) and (2) are inverse to each other.
\end{proof}
\begin{Rem}
  Using \eqref{epsepst} in \nmref{keyrem}, we can describe
  the tensor product of two left $H$-comodules $M$ and $N$, 
  which is their tensor product $M\ou RN$, by the isomorphic
  subspace 
  $$\{\epsilon(m\swm 1n\swm 1)m\sw 0\o n\sw 0|m\in M,n\in N\}\subset M\o N$$
  since we have
  \begin{multline*}
    me\so 1\o e\so 2n
      =\epsilon(m\swm 1\epsilon_t(1\sw 1))\epsilon(1\sw 2n\swm 1)\o m\sw 0\o n\sw 0
      \\=\epsilon(m\swm 11\sw 1)\epsilon(1\sw 2n\swm 1)m\sw 0\o n\sw 0
      =\epsilon(m\swm 1n\swm 1)m\sw 0\o n\sw 0.
  \end{multline*}
  for $m\in M$ and $n\in N$. Both versions of a tensor product
  in the category of $H$-comodules (the tensor product over $R$ and
  the version that is a subspace of the tensor product over $k$)
  were discussed and compared by Nill \cite[Sec.4]{Nil:AWB} in the
  case where $H$ is finite dimensional.
\end{Rem}

\end{document}